\documentclass[11pt]{amsart}
\usepackage{mathrsfs,url,amssymb}

\theoremstyle{plain}
    \newtheorem{thm}{Theorem}
    
    \newtheorem{lem}[thm]{Lemma}
    \newtheorem{prop}[thm]{Proposition}

\theoremstyle{definition}
    \newtheorem{defn}[thm]{Definition}
    
\theoremstyle{remark}
    \newtheorem{rem}[thm]{Remark}
\theoremstyle{remark}

\newcommand{\To}{\rightarrow}

\newcommand{\nin}{\notin}
\newcommand{\sm}{\setminus}
\renewcommand{\c}{X \setminus}
\newcommand{\inv}{^{-1}}

\newcommand{\cl}[1]{\langle #1 \rangle}

\newcommand{\cls}[1]{\cl{\{#1\}\cup\S}}

\renewcommand{\O}{{\mathscr O}}
\newcommand{\On}{{\mathscr O}^{(n)}}

\newcommand{\Oo}{{\mathscr O}^{(1)}}

\newcommand{\Pow}{{\mathscr P}}
\renewcommand{\P}{{\mathscr P}}

\newcommand{\C}{{\mathscr C}}
\newcommand{\R}{{\mathscr R}}
\newcommand{\F}{{\mathscr F}}

\newcommand{\I}{{\mathscr I}}
\newcommand{\A}{{\mathscr A}}
\newcommand{\B}{{\mathscr B}}
\newcommand{\D}{{\mathscr D}}
\newcommand{\K}{{\mathscr K}}
\newcommand{\G}{{\mathscr G}}
\renewcommand{\S}{{\mathscr S}}
\newcommand{\M}{{\mathscr M}}

\newcommand{\rest}{\upharpoonright}

\author[M.\,Pinsker]{Michael Pinsker}
\address{Algebra\\TU Wien\\Wiedner Hauptstra\ss e 8-10/104\\A-1040 Wien, Austria}
\email{marula@gmx.at} \urladdr{http://www.dmg.tuwien.ac.at}
\title[Unary clones containing the permutations]{The number of unary clones containing the permutations on an infinite set}


\thanks{Support by DOC [Doctoral Scholarship Programme of the Austrian
Academy of Sciences], and later by the Postdoctoral Fellowship of
the Japan Society for the Promotion of Science (JSPS) is
gratefully acknowledged.}

\begin{document}

\begin{abstract}
    We calculate the number of unary clones (submonoids of the full transformation
    monoid) containing the permutations, on an infinite base set. It turns out
    that this
    number is quite large, on some cardinals as large as the whole clone
    lattice. Moreover we find that, with one exception, even the cardinalities of the intervals
    between the monoid of all permutations and the maximal submonoids
    of the full transformation monoid are as large. Whether or not the only exception is of the same cardinality
    as the other intervals depends on additional axioms of set theory.
\end{abstract}

\maketitle

\section{Background and the result}
    Fix a set $X$ and consider for all $n \geq 1$ the set $\On$ of
    $n$-ary operations on $X$. If we take the union $\O=\bigcup_{n\geq
    1}\On$ over these sets, we obtain the set of all operations on $X$ of finite arity. A
    \emph{clone} is a subset of $\O$ which contains all functions
    of the form $\pi^n_k(x_1,\ldots,x_n)=x_k$ ($1\leq k\leq n$), called the projections, and
    which is closed under composition of functions. With the order of set-theoretical inclusion,
    the clones on $X$ form a
    complete algebraic lattice $Cl(X)$. We wish to describe this lattice
    for infinite $X$, in which case it has cardinality
    $2^{2^{|X|}}$.

    A clone is called \emph{unary} iff it contains only essentially unary functions, i.e., functions
    which depend on only one variable. Unary clones correspond in an obvious way to submonoids of the
    full transformation monoid $\Oo$ and we shall not distinguish between the two notions in the
    following. We say that a unary clone $\C\neq \Oo$ is
    \emph{precomplete} or \emph{maximal}
    iff $\C$
    together with any unary function $f\in\Oo\setminus\C$
    generates
    $\Oo$, i.e. iff the smallest clone containing $\C$ as well as
    $f$ is $\Oo$.
    In \cite{Pin033}, the author determined all precomplete submonoids of the
    full transformation monoid $\Oo$ that contain the
    permutations for all infinite $X$, which was a generalization
    from the countable (\cite{Gav65}). The number of such clones turned out to be rather
    small compared with the size of the clone lattice:
    On an infinite set $X$ of cardinality $\aleph_\alpha$ there exist
    $2|\alpha|+5$ precomplete unary clones, so in particular there are only five precomplete unary clones on
    a countably infinite set $X$.

    \begin{thm}\label{THM:allPrecompleteMonoids}
        Let $X$ be an infinite set of cardinality $\kappa$. If $\kappa$ is regular, then the
        precomplete submonoids of $\Oo$ that contain the permutations
        are exactly the monoid $\A$ and the monoids $\G_\xi$ and $\M_\xi$ for
        $\xi=1$ and $\aleph_0\leq \xi\leq\kappa$, $\xi$ a cardinal, where
        \begin{itemize}
            \item{$\A=\{f\in\Oo:f^{-1}[\{y\}]$ is small for
            almost all $y\in X\}$}
            \item{$\G_\xi=\{f\in\Oo: f\text{ is
            }\xi\text{-injective or not
            }\xi\text{-surjective}\}$}
            \item{$\M_\xi=\{f\in\Oo: f \text{ is
            }\xi\text{-surjective or not
            }\xi\text{-injective}\}$}
        \end{itemize}
        If $\kappa$ is singular, then the same is true
        with the monoid $\A$ replaced by
        \begin{itemize}
            \item $\A'=\{f\in\Oo: \exists \xi < \kappa \,\,(\,|f\inv[\{x\}]|\leq\xi\text{ for almost all }x\in
            X\,)\,\}.$
        \end{itemize}
    \end{thm}

        In the theorem, a set is \emph{small} iff it has cardinality smaller
    than the cardinality of $X$, a property holds for \emph{almost all} $y\in X$
    iff it holds for all $y\in X$ except for a small set, a
    function $f\in\Oo$ is \emph{$\xi$-surjective} iff $|X\sm
    f[X]|<\xi$, and it is \emph{$\xi$-injective} iff there is a set
    $Y\subseteq X$ of cardinality smaller than $\xi$ such that the
    restriction of $f$ to $X\sm Y$ is injective.

    With this result, the
    question arose whether it was possible to describe the whole
    interval $[\S,\Oo]$ of the clone lattice,
    where $\S$ is the set of permutations
    of $X$. We
    show that compared to the number of its dual atoms, this
    interval is quite large. In particular,
    on a countably infinite set $X$ it equals the size of the
    whole clone lattice.

    \begin{thm}\label{THM:numberOfMonoidsAboveS}
        Let $X$ be an infinite set of cardinality
        $\kappa=\aleph_\alpha$. Then there exist $2^{2^\lambda}$ submonoids of $\Oo$ which contain
        all permutations, where
        $\lambda=\max\{\,|\alpha|,\aleph_0\}$. Moreover, if $\kappa$
        is regular, then $|[\S,\G]|=2^{2^\lambda}$ for every precomplete monoid above
        $\S$; in fact, $|[\S,\D]|=2^{2^\lambda}$, where $\D$ is the
        intersection of the precomplete elements of $[\S,\Oo]$. If
        $\kappa$ is singular, then $|[\S,\G]|=2^{2^\lambda}$ for all precomplete monoids
        except $\A'$: If $\lambda<\kappa$, then $|[\S,\A']|=|[\S,\D]|=2^{2^\lambda}$,
        but if $\lambda=\kappa$, then $|[\S,\A']|=|[\S,\D]|=2^{(\kappa^{<\kappa})}$ (where $\kappa^{<\kappa}=\sup\{\kappa^\xi:\xi<\kappa\}$).
    \end{thm}

    \subsection{Notation}
        For any set $Y$, we denote the power set of $Y$ by
        $\P(Y)$. The smallest clone containing a set of functions
        $\F\subseteq\O$ is denoted by $\cl{\F}$.
        If $f\in\Oo$, we write $ker(f)\subseteq\P(X)$ for the kernel of
        $f$.
\section{The proof of Theorem \ref{THM:numberOfMonoidsAboveS}}
    \begin{defn}
        Set $\K=\{\xi:\xi \text{ a cardinal and }
        \xi\leq\kappa\}$; then $|\K|=\lambda$. Define for every $f\in\Oo$ a function
        $$
            s_f:\quad\begin{matrix} \K &\To& \K\\
            \xi&\mapsto& |\{A\in ker(f): |A|=\xi\}| \end{matrix}
        $$
        In words, the function assigns to every $\xi\leq\kappa$ the
        number of equivalence classes in the kernel of $f$ which have cardinality $\xi$. We call $s_f$
        the \emph{kernel sequence} of $f$.
    \end{defn}
    \begin{lem}
        If $f,g\in\Oo$ are unary functions satisfying $s_f=s_g$
        and $|\c f[X]|=|\c g[X]|$, then there exist $\beta,
        \gamma\in\S$ such that $f=\beta\circ g\circ\gamma$.
    \end{lem}
    \begin{proof}
        The assumption $s_f=s_g$ implies that there is $\gamma\in\S$ such that $ker(f)=ker(g\circ\gamma)$.
        Obviously, $|f[X]|=|g[X]|=|g\circ\gamma[X]|$ as $s_f=s_g$.
        Together with the fact that $|\c f[X]|=|\c g[X]|$ this
        implies that we can find $\beta\in\S$ such that
        $f[X]=\beta\circ g\circ\gamma[X]$, and since $ker(f)=ker(g\circ\gamma)$ also so that
        $f=\beta\circ g\circ\gamma$.
    \end{proof}
    \begin{prop}
        The number $\mu$ of submonoids of $\Oo$ containing $\S$ is at
        most $2^{2^\lambda}$.
    \end{prop}
    \begin{proof}
        By the preceding lemma, the clone a function $f\in\Oo$ generates together with $\S$
        is determined by $s_f$ and the cardinality of $X\sm f[X]$. There exist at most
        $\lambda^\lambda$ different kernel
        sequences and $\lambda$ possibilities for the cardinality of the complement of the range of a function in
        $\Oo$. Thus, modulo $\S$ there are only
        $\lambda^\lambda \cdot \lambda=\lambda^\lambda=2^\lambda$
        different functions in $\Oo$. Therefore, $\mu\leq 2^{2^\lambda}$.
    \end{proof}

    We will now show the other inequality. Fix any sequence
    $(n_i)_{i\in\omega}$ of natural numbers such that $\sum_{j<i}
    n_j<n_i$ for all $i\in\omega$. Set
    $\R=\{n_i\}_{i\in\omega}\cup\{\xi\in\K:\xi\text{ infinite
    successor}\}$. Then $|\R|=|\K|=\lambda$. For all $f\in\Oo$,
    write $\widetilde{s}_f=s_f\rest_\R$ for the restriction of its
    kernel sequence to $\R$.

    Observe that for all $\xi\in\R$ we have that
    $\sum_{\eta<\xi,\eta\in\R}\eta<\xi$: For $\xi$ finite, this is
    because we chose the finite elements of $\R$ that way, and if
    $\xi$ is infinite, then it is a successor cardinal so that the
    left side of the inequality is clearly bounded by its predecessor.

    We say that $A\subseteq \R$ is  \emph{unbounded} iff $\sum_{\xi\in A}\xi=\kappa$. Assign
    to every unbounded $A\subseteq\R$ a function $f_A\in\Oo$ satisfying
    $s_{f_A}(\xi)=1$ whenever $\xi\in A$, and $s_{f_A}(\xi)=0$ whenever $\xi\in\K\sm
    A$. The fact that $A$ is unbounded guarantees the existence of
    $f_A$.
    \begin{lem}
        If $A\subseteq\R$ is unbounded and $g\in \Oo$, then $\widetilde{s}_{g\circ
        f_A}\leq \widetilde{s}_{f_A}$.
    \end{lem}
    \begin{proof}
        Consider an arbitrary $B\in ker(g\circ
        f_A)$ with $|B|=\xi\in\R$. We claim there exists $C\subseteq B$ of
        cardinality $\xi$ such that $C\in ker(f_A)$. For suppose
        to the contrary this is not the case. Being an element of
        $ker(g\circ f_A)$, $B$ is the union of sets in the kernel
        of $f_A$: $B=\bigcup_{i\in\delta} B_i$, for $B_i\in ker(f_A)$ and some ordinal
        $\delta$. By our assumption, $|B_i|<\xi$ for all
        $i\in\delta$. Thus, $|B|=|\bigcup_{i<\delta} B_i|\leq
        \sum_{D\in ker(f_A), |D|<\xi}|D|=\sum_{\eta\in
        A,\eta<\xi}\eta<\xi$, contradiction. So for all $B\in ker(g\circ
        f_A)$ with $|B|=\xi\in\R$ we injectively find $C\in ker(f_A)$ of the same
        cardinality, which proves the lemma.
    \end{proof}
    \begin{lem}
        Let $A,A_1,\ldots,A_n\subseteq\R$ be unbounded and such that $A\nsubseteq A_i$ for all $1\leq i\leq
        n$. Then $f_A\nin\cls{f_{A_1},\ldots,f_{A_n}}$.
    \end{lem}
    \begin{proof}
        Clearly, every unary $t\in\cls{f_{A_1},\ldots,f_{A_n}}$ which is not a permutation has a representation of the
        form $t=g\circ f_{A_i}\circ\beta$, where $g\in\Oo$,
        $\beta\in\S$ and $1\leq i\leq n$. But then
        $\widetilde{s}_t\leq \widetilde{s}_{f_{A_i}}$ by the preceding lemma, so that
        $\widetilde{s}_t\neq
        \widetilde{s}_{f_A}$ and therefore $t\neq f_A$.
    \end{proof}

    It is a fact that if $Y$ is any set, then there
    exists a family
    $\I\subseteq\P(Y)$ such that $|\I|=|\P(Y)|=2^{|Y|}$ and such that the sets of
    $\I$ are pairwise \emph{incomparable},
    i.e., $A\nsubseteq B$ holds for all distinct $A, B\in\I$. For example, it is
    a well-known theorem of Hausdorff that
    there exist \emph{independent} families of subsets of $Y$
    of that size,
    where $\I\subseteq\P(Y)$ is called independent iff every nontrivial Boolean combination of sets from $\I$ is
    nonempty,
    i.e., whenever $\B_1,\B_2\subseteq\I$ are finite, nonempty and disjoint, then
    $$
        \bigcap_{A\in\B_1}A\,\,\cap\bigcap_{A\in\B_2} (Y\sm
        A)\neq\emptyset.
    $$
    See the textbook \cite[Lemma 7.7]{Jec02} for a proof of this.

    There is an independent family of unbounded subsets of $\R$
    which has cardinality ${2^\lambda}$: If $\I\subseteq\P(\R)$ is
    independent of size $2^\lambda$, then either $\I$ or $\I'=\{\R\sm
    A: A\in \I\}$ contains $2^\lambda$ unbounded sets, the family of which is independent.

    \begin{prop}\label{PROP:atLeast22lambda}
        There is an order embedding from $\P(2^\lambda)$ into $[\S,\Oo]$. In particular,
        the number $\mu$ of submonoids of $\Oo$ containing $\S$ is at
        least $2^{2^\lambda}$.
    \end{prop}
    \begin{proof}
        Let $\I\subseteq\P(\R)$ be an independent family
        of unbounded subsets of $\R$ with $|\I|=2^\lambda$.
        Define for every $\B\subseteq\I$ a monoid
        $\C_\B=\cls{f_A:A\in\B}$. Then for all $\B_1,\B_2\subseteq\I$
        we have that if $\B_1\nsubseteq\B_2$, then
        $\C_{\B_1}\nsubseteq \C_{\B_2}$: Indeed, by the preceding lemma $f_A\in\C_{\B_1}\sm\C_{\B_2}$ for any
        $A\in\B_1\sm\B_2$. Together with the fact that larger subsets of $\I$ yield larger
        clones, this implies that
        the mapping $\varphi: \Pow(\I)\To [\S,\Oo]$ assigning to every $\B\subseteq\I$ the clone
        $\C_\B$ is an order embedding.
        Hence, there exist $|\Pow(\I)|=2^{2^\lambda}$ distinct
        monoids containing the permutations.
    \end{proof}
        This completes the proof of the first statement of Theorem
        \ref{THM:numberOfMonoidsAboveS}.
    \begin{prop}\label{PROP:additionalStatementReg}
        Let $\kappa$ be regular and let $\D$ be the intersection of the precomplete submonoids of $\Oo$ containing
        $\S$. There is an order embedding of $\P(2^\lambda)$ into $[\S,\D]$. Hence, $|[\S,\D]|=2^{2^\lambda}$.
    \end{prop}
    \begin{proof}
        Since in the proof of Proposition
        \ref{PROP:atLeast22lambda} we considered only functions
        $f\in\Oo$ with $s_f(\kappa)\leq 1$, all those functions
        were elements of $\A$. Also, we did not care about the
        size of the complement of the range of $f$; if we assume it to be of
        cardinality $\kappa$, then all functions of the
        construction are not $\kappa$-surjective and therefore elements of $\G_\xi$, for all cardinals
        $\xi=1$ and $\aleph_0\leq\xi \leq\kappa$. Since for any unbounded $A\subseteq \R$ and any
        small $Y\subseteq X$ there is $\xi\in A$ with $\xi> |Y|$, the
        fact that $f_A$ has a class of size $\xi$ in its kernel
        yields that $f_A$ is not injective on the complement of $Y$.
        Therefore, the $f_A$ used in the construction are not
        $\kappa$-injective and hence are elements of $\M_\xi$, for all
        $\xi=1$ and $\aleph_0\leq\xi\leq\kappa$. This proves the proposition.
    \end{proof}
        We now turn to the case when $\kappa$ is singular. The argument
        of the preceding proposition yields
    \begin{prop}\label{PROP:additionalStatementSing}
        Let $\kappa$ be singular and let $\G\neq\A'$ be a precomplete submonoid of $\Oo$ containing
        $\S$. There is an order embedding of $\P(2^\lambda)$ into $[\S,\G]$.
        In particular, $|[\S,\G]|=2^{2^\lambda}$.
    \end{prop}

    \begin{prop}\label{PROP:singularEasyCase}
        Let $\kappa$ be singular such that $\lambda<\kappa$.
        Then $|[\S,\A']|=2^{2^\lambda}$.
    \end{prop}
    \begin{proof}
        Since the functions used
        in our construction satisfy $s_f(\xi)\leq 1$ for all
        $\xi\in\K$, we have $|f[X]|\leq\lambda<\kappa$, and hence
        $|f\inv[\{x\}]|=0$ for almost all $x\in X$; therefore those functions are elements of
        $\A'$. Hence, $|[\S,\A']|=2^{2^\lambda}$.\\
    \end{proof}
    \begin{prop}\label{PROP:singularConplicatedCase}
        Let $\kappa$ be singular such that $\lambda=\kappa$. Then $|[\S,\A']|=2^{(\kappa^{<\kappa})}$.
    \end{prop}
    \begin{proof}
        We first calculate the number of different kernel sequences of functions in
        $\A'$. Let $s_f:\K\To \K$ be such a sequence; then $f\in\A'$
        iff there is $\xi<\kappa$ such that
        $\sum_{\xi\leq\eta\leq\kappa}s_f(\eta)=\tau<\kappa$.
        Fixing $\xi$ and $\tau$, we have $\kappa^\tau$ possibilities for
        the part of $s_f$ between $\xi$ and $\kappa$. Taking the sum over all $\tau<\kappa$,
        we obtain $\kappa^{<\kappa}$ possibilities for $s_f$ between $\xi$ and $\kappa$.
        Since below $\xi$ there are no conditions on $s_f$ in order to make $f$ an element of $\A'$,
        there are exactly $\kappa^\xi$ possibilities for the restriction of $s_f$ to $\xi$,
        so that we have a total of
        $\kappa^\xi+\kappa^{<\kappa}$ kernel sequences of
        functions $f$ with $\sum_{\xi\leq\eta\leq\kappa}s_f(\eta)<\kappa$.
        Since $\xi<\kappa$ can be arbitrary, we take the sum over all $\xi<\kappa$ and find
        that there are $\kappa^{<\kappa}$ distinct kernel sequences of
        functions in $\A'$. Hence, $|[\S,\A']|\leq 2^{(\kappa^{<\kappa})}$.

        \textbf{Claim.} There exists a family $\I$ of pairwise incomparable small
        unbounded subsets of $\R$ which has cardinality $\kappa^{<\kappa}$.

        To prove this, we first observe that for all $cf(\kappa)\leq \xi<\kappa$ there
        exists a family $\I_\xi$ of pairwise incomparable unbounded subsets of $\R$ of
        cardinality $\xi$ such that $|\I_\xi|=\kappa^\xi$ ($cf(\kappa)$ denotes the cofinality of $\kappa$).
        Indeed,
        write $\R=\R'\cup\R''$, where $\R'$ and $\R''$ are disjoint,
        and $\R''$ is unbounded and of cardinality $\xi$.
        Now let $\I_\xi'$ be a family of pairwise
        incomparable subsets
        of $\R'$ of cardinality $\xi$ with $|\I_\xi'|=\kappa^\xi$.
        To see that $\I_\xi'$ exists, observe
        that every function $f\in\kappa^\xi$ is a subset of
        $\xi\times\kappa$, and that all those functions are
        incomparable as subsets of $\xi\times\kappa$. Thus a family of size
        $\kappa^\xi$ of
        pairwise incomparable sets of size $\xi$ exists on $\xi\times\kappa$,
        and therefore also on $\R'$ since $|\R'|=|\xi\times\kappa|=\kappa$.
        Now we set $\I_\xi=\{A\cup \R'': A\in \I_\xi'\}$ to obtain
        the
        family $\I_\xi$ having the desired properties. Finally to prove the claim, write $\R$ as a disjoint union
        $\R=\bigcup_{cf(\kappa)\leq \xi< \kappa}\R_\xi$
        of sets $\R_\xi$ of cardinality $\kappa$ (which also implies that they are unbounded). Fix a family
        $\I_\xi$ of pairwise incomparable unbounded subsets of $\R_\xi$ of
        cardinality $\xi$ such that
        $|\I_\xi|=\kappa^\xi$, for all
        $\xi$. Then the family $\I=\bigcup_{cf(\kappa)\leq \xi < \kappa}  \I_\xi$
        consists of pairwise incomparable small unbounded subsets
        of $\R$ and has cardinality $\kappa^{<\kappa}$.

        Having small range, the functions corresponding to the sets in $\I$ are all members of $\A'$,
        so that we obtain
        $2^{(\kappa^{<\kappa})}$ clones in the interval
        $[\S,\A']$.

    \end{proof}
    \begin{prop}\label{PROP:finalStatementSing}
        Let $\kappa$ be singular and let $\D$
        be the intersection of the precomplete submonoids of $\Oo$ containing
        $\S$. If $\lambda<\kappa$, then $|[\S,\D]|=2^{2^{\lambda}}$.
        If $\lambda=\kappa$, then $|[\S,\D]|=2^{(\kappa^{<\kappa})}$.
    \end{prop}
    \begin{proof}
        One only needs to combine the proofs of Propositions
        \ref{PROP:additionalStatementReg},
        \ref{PROP:additionalStatementSing},
        \ref{PROP:singularEasyCase},
        and
        \ref{PROP:singularConplicatedCase}; we
        leave the details to the reader.
    \end{proof}
    \begin{rem}
        If GCH holds, then $2^{(\kappa^{<\kappa})}=2^{2^\kappa}$, so
        in this case we have $[\S,\D]=2^{2^\lambda}$ on all
        infinite $X$. However, for any singular $\kappa$ it is is also consistent that $2^\kappa<
        2^{(\kappa^{<\kappa})}<2^{2^{\kappa}}$. Therefore, if $\kappa$ is singular and $\lambda=\kappa$,
        then the intervals $[\S,\A']$ and $[\S,\D]$ can be smaller than
        $2^{2^{\lambda}}$. In particular we have that whether or not the intervals $[\S,\A']$
        and, say, $[\S,\M_1]$ are of equal cardinality depends on
        the set-theoretical universe.
    \end{rem}

\bibliographystyle{alpha}

\end{document}